# A new family of *q*-Bernstein polynomials: Probabilistic viewpoint

**Ayse Karagenc[1], Mehmet Acikgoz[2] and Serkan Araci[3,*]**

[1,2]Department of Mathematics, Faculty of Science and Arts, Gaziantep University, TR-27310 Gaziantep, Türkiye

[3]Department of Computer Engineering, Faculty of Engineering, Hasan Kalyoncu University, TR-27010 Gaziantep, Türkiye

E-Mails: ayse_karagenc@hotmail.com[1]; acikgoz@gantep.edu.tr[2]; serkan.araci@hku.edu.tr[3]

*Corresponding Author

**Abstract.** In this paper, we introduce a new class of polynomials, called probabilistic *q*-Bernstein polynomials, alongside their generating function. Assuming *Y* is a random variable satisfying moment conditions, we use the generating function of these polynomials to establish new relations. These include connections to probabilistic Stirling numbers of the second kind and higher-order probabilistic Bernoulli polynomials associated with *Y*. Additionally, we derive recurrence and differentiation properties for probabilistic *q*-Bernstein polynomials. Utilizing Leibniz's formula, we give an identity for the generating function of these polynomials. In the latter part of the paper, we explore applications by choosing appropriate random variables such as Poisson, Bernoulli, Binomial, Geometric, Negative Binomial, and Uniform distributions. This allows us to derive relationships among probabilistic *q*-Bernstein polynomials, Bell polynomials, Stirling numbers of the second kind, higher-order Frobenius-Euler numbers, and higher-order Bernoulli polynomials. We also present *p*-adic *q*-integral and fermionic *p*-adic *q*-integral representations for probabilistic *q*-Bernstein polynomials.

1. Introduction

    1.1 Introduction of Probabilistic viewpoint

Assume that *Y* is a random variable satisfying appropriate moment conditions. Adell and Lekuona [3] introduced a different generalization of Stirling numbers of the second kind as

$$\sum_{n=m}^{\infty} S_Y(n,m) \frac{t^n}{n!} = \frac{\left(\mathrm{E}\left[e^{tY}\right]-1\right)^m}{m!},$$

which is also called probabilistic Stirling numbers of the second kind. Note that alternating notation for probabilistic Stirling numbers of the second kind is $S_Y(n,m) := \begin{Bmatrix} n \\ m \end{Bmatrix}_Y$. In the case $Y = 1$, one has

$$\sum_{n=m}^{\infty} S(n,m) \frac{t^n}{n!} = \frac{(e^t - 1)^m}{m!},$$

where $S(n,m)$ is called Stirling numbers of the second kind. Obviously that $S_{Y=1}(n,m) := \begin{Bmatrix} n \\ m \end{Bmatrix}_{Y=1} := S(n,m) := \begin{Bmatrix} n \\ m \end{Bmatrix}$. Thus, various probabilistic generalizations of special functions and polynomials have been studied by several mathematicians, as detailed in the following section.

### 1.2 Literature Review

In the year 2019, Adell and Lekuona [2] initiated the different generalization of Stirling numbers of the second kind, called probabilistic Stirling numbers of the second kind. Later, Adell [3] gave a new generalization of the Stirling numbers of the second kind associated with each complex-valued random variable satisfying appropriate integrability conditions, and derived the applications such as determining asymptotic behavior without utilizing central limit theorem, and Lévy processes and cumulants. After that, Kim proceeded the probabilistic viewpoint to several special numbers and polynomials such as probabilistic degenerate Bell polynomials associated with random variables [14], probabilistic Bernoulli and Euler polynomials [15], probabilistic degenerate Stirling polynomials of the second kind [16], probabilistic type 2 Bernoulli and Euler polynomials [21], probabilistic central Bell polynomials [30] and probabilistic degenerate Fubini polynomials associated with random variables [31]. Also, recent investigations including probabilistic properties of polynomials defined over the ring of $p$-adic integers under the Haar measure [19], Appell polynomials subject to a random variable [27], a probabilistic approach in order to compute the Adomian polynomials [29] and various applications of several special functions and polynomials including Hermite based Genocchi polynomials [33] and generalized Bell-Based Appell polynomials [34] have been studied extensively.

### 1.3 Plan of this paper

Let $Y$ be a random variable fulfilling the moment criteria by

$$E[|Y|^n] < \infty \quad (n \in \mathbb{N}_0) \quad \text{and} \quad \lim_{n \to \infty} \frac{|v|^n E[|Y|^n]}{n!} = 0, \quad (|v| < r; r > 0),$$

where E means mathematical expectation. Henceforth, we see (*cf.* [3]) that

$$E[e^{ivY}] = \sum_{n=0}^{\infty} E[Y^n] i^n \frac{v^n}{n!}, \quad (i = \sqrt{-1}; |v| < r).$$

Alternatively, (see [3])

$$E[e^{v|Y|}] < \infty, \quad |v| < r.$$

Let $\{Y_j\}_{j=1}^{l} = (Y_1, Y_2, \ldots, Y_l)$ be a sequence of mutually independent samples of $Y$ and be denoted by

$$S_l = Y_1 + Y_2 + \cdots + Y_l, \quad (l \in \mathbb{N}).$$

with the initial assumption $S_0 := 0$.

The six distributions used in this paper are listed below (see reference [28]), assuming the moment generating functions are written within their appropriate intervals of convergence:

1. ***Poisson distribution.*** $Y \sim Poisson(\alpha)$ with the parameter $\alpha > 0$ with the probability density function (pdf) $P(Y = y) = f(y) = \frac{e^{-\alpha} \alpha^y}{y!}, y \in \mathbb{N}_0$ yielding moment generating function (mgf) as

$$E[e^{Yv}] = e^{\alpha(e^v - 1)}. \tag{1.1}$$

2. ***Bernoulli distribution.*** $Y \sim Bernoulli(p_1)$ with the success probability $p_1$ with pdf $f(y) = p_1^y (1 - p_1)^{1-y}, (y \in \{0,1\})$ yielding mgf

$$E[e^{Yv}] = 1 - p_1 + p_1 e^v. \tag{1.2}$$

3. ***Binomial distribution.*** $Y \sim b(n, p_1)$ with pdf $f(y) = \binom{n}{y} p_1^y (1 - p_1)^{n-y}, (y = 0,1,2, \ldots, n)$ yielding mgf

$$E[e^{Yv}] = (1 - p_1 + p_1 e^v)^n. \tag{1.3}$$

4. ***Geometric distribution.*** $Y \sim Geo(p_1)$ with pdf $f(y) = p_1^y (1 - p_1)^{y-1}, (y = 1,2, \ldots, n)$ yielding mgf

$$E[e^{Yv}] = \frac{p_1 e^v}{1 - (1 - p_1) e^v}. \tag{1.4}$$

5. ***Negative Binomial distribution.*** $Y \sim NB(n, p_1)$ with pdf

$$f(y, a, p_1) = \binom{y-1}{a-1} p_1^a (1 - p_1)^{y-a}, \quad y = a, a+1, \ldots \tag{1.5}$$

yielding mgf

$$E[e^{Yv}] = \left(\frac{p_1 e^v}{1-(1-p_1)e^v}\right)^a \tag{1.6}$$

6. **Uniform distribution**. Let $Y \sim U(0,1)$ be a random variable (continuous random variable) with pdf

$$f(y) = 1, \quad 0 \leq y \leq 1$$

$$E[e^{Yv}] = \frac{e^v - 1}{v}. \tag{1.7}$$

Recently, in (see [8]), the authors introduced a probabilistic version of classical Bernstein polynomials, denoted by $B_{k,n}^Y(x)$, using the following generating function:

$$\sum_{n=k}^{\infty} B_{k,n}^Y(x) \frac{t^n}{n!} = \frac{(tx)^k}{k!} (E[e^{tY}])^{1-x}, \tag{1.8}$$

and derived recurrence, derivative and relation formulas linking these polynomials to other special polynomials such as Euler polynomials, Bernoulli polynomials of higher order, Frobenius-Euler polynomials of higher order, Stirling numbers of the second kind and Bell polynomials. When $Y = 1$ in Eq. (1.8), we obtain

$$\sum_{n=k}^{\infty} B_{k,n}(x) \frac{t^n}{n!} = \frac{(tx)^k}{k!} e^{(1-x)t},$$

which corresponds to the generating function for the Bernstein polynomials. This leads us to pose the following question:

*Is it possible to construct a further generalization of the probabilistic Bernstein polynomials within the framework of q-analysis?*

We first consider $q$-Bernstein polynomials, which were also introduced by Simsek and Acikgoz [24, Eq. (3.8)]:

$$\sum_{n=k}^{\infty} Y_n(k; x; q) \frac{v^n}{n!} = \frac{(v[x]_q)^k}{k!} e^{[1-x]_q v},$$

which was also studied, and different proof was provided by Kim *et al.* [18, Eq. (2.1)], by taking $Y_n(k;x;q) = B_{k,n}(x,q)$. Throughout this paper, we will use $B_{k,n}(x,q)$, as it closely resembles the classical one. Here, $[x]_q$ is $q$-extension of $x$, defined as

$$[x]_q := \frac{q^x - 1}{q - 1}, \quad (0 < q < 1).$$

Note that $\lim_{q \to 1^-} [x]_q = x$.

Our aim is to give an answer to the question above by introducing the following generating function of the probabilistic $q$-Bernstein polynomials:

$$F_{r;Y}(x,v,q) = \sum_{n=r}^{\infty} B_{r,n}^Y(x,q) \frac{v^n}{n!} = \frac{\left(v[x]_q\right)^r}{r!} \left(\mathrm{E}\left[e^{vY}\right]\right)^{[1-x]_q}, \quad \left(0 \le x \le 1; r \in \{0,1,2,\cdots,n\}\right), \qquad (1.9)$$

where $B_{r,n}^Y(x,q)$ yields $B_{r,n}(x,q)$ in the case when $Y = 1$, that is, $B_{r,n}^{Y=1}(x,q) := B_{r,n}(x,q)$. By introducing a new and different family of Bernstein polynomials as given in Eq. (1.9), we will derive some new interesting identities, relations and results. Also, by considering appropriate random variable $Y$ such as Bernoulli, Poisson, Binom, Geometric, Negative Binom, and Uniform distributions, we explored interesting results among some special polynomials such as Bell polynomials, Bernoulli polynomials of higher order, Frobenius-Euler polynomials of higher order, Stirling numbers of the second kind, etc., which we stated in Section 3.

**1.4 Notations**

Here we briefly summarize some mathematical tools which will be used in deriving main results of this paper. We now begin with the following notations:

$$\mathbb{N} := \{1,2,3,\ldots\}, \mathbb{N}_0 := \{0,1,2,3,\ldots\} = \mathbb{N} \cup \{0\}.$$

$\mathbb{Z}, \mathbb{Q}, \mathbb{R}$, and $\mathbb{C}$ will respectively, be denoted integers, rational numbers, real numbers and complex numbers.

For $0 < q < 1$, we recall $q$-extension of $x$, denoted by $[x]_q$, as follows:

$$[x]_q = \frac{q^x - 1}{q - 1}.$$

Note that $\lim_{q \to 1^-} [x]_q = x$. Additionally, we list the following properties which will be shown easily by making use of the definition $q$-extension of $x$.

1. $[x - y]_q = [x]_q - q^{x-y}[y]_q$,
2. $[-x]_q = -q^{-x}[x]_q$, $[x]_{q^{-1}} = q^{1-x}[x]_q$,
3. $[1 - x]_q = 1 - [x]_{q^{-1}}$.

Then we have

$$[1-x]_q^m = q^{-mx} \sum_{l=0}^{m} \binom{m}{l} (-1)^l [x]_q^l. \tag{1.10}$$

See [9], [11], [12], [32] for more details on $q$-analysis

Let $B_{r,n}(x)$ represent the Bernstein polynomials of degree $n$ (also known as classical Bernstein polynomials, Bernstein polynomials, or Bernstein polynomial bases) defined by the following generating function. For $x \in [0,1]$, $n \in \mathbb{N}_0 := \mathbb{N} \cup \{0\}$ and $v \in \mathbb{C}$, we have

$$F_r(x,v) = \frac{(vx)^r}{r!} e^{v(1-x)} = \sum_{n=r}^{\infty} B_{r,n}(x) \frac{v^n}{n!}, \tag{1.11}$$

where

$$B_{r,n}(x) = \binom{n}{r} x^r (1-x)^{n-r}, \qquad r = 0, 1, \dots, n.$$

For further details on Bernstein polynomials, please refer to [1, 24].

Each Bernoulli polynomial of degree $n$, $B_n(x)$, can be determined using the generating function:

$$\frac{v}{e^v - 1} e^{xv} = \sum_{n=0}^{\infty} B_n(x) \frac{v^n}{n!}, \qquad (|v| < 2\pi), \tag{1.12}$$

where, with $x = 0$ in Eq. (1.12), $B_n(0) := B_n$, representing the Bernoulli numbers (see [5,6,13,15,26]).

In [13], the Bernoulli polynomials of higher order, $B_n^{(\alpha)}(x)$, can be computed by

$$\left(\frac{v}{e^v - 1}\right)^\alpha e^{xv} = \sum_{n=0}^{\infty} B_n^{(\alpha)}(x) \frac{v^n}{n!}, \qquad (|v| < 2\pi). \tag{1.13}$$

The Euler polynomials, $E_n(x)$, described in In [5,6,10,12,15,26], are given by

$$\frac{2}{e^v + 1} e^{xv} = \sum_{n=0}^{\infty} E_n(x) \frac{v^n}{n!}, \qquad (|v| < \pi), \tag{1.14}$$

where, for $x = 0$ in (1.14), $E_n(0) := E_n$ denotes Euler numbers.

In [14], Bell polynomials, $\phi_n(x)$, are defined by

$$\sum_{n=0}^{\infty} \phi_n(x) \frac{v^n}{n!} = e^{x(e^v - 1)},$$

which corresponds to the moment generating function of the Poisson distribution with mean $x$.

Lastly, the Frobenius-Euler polynomials of order $r \in \mathbb{R}$, as described in in [4], are given by

$$\left(\frac{1-u}{e^v - u}\right)^r e^{xv} = \sum_{n=0}^{\infty} H_n^{(r)}(x|u) \frac{v^n}{n!}, \qquad (u \in \mathbb{C} - \{1\}).$$

We now review several $p$-adic tools and notations that will be used in the final section of this paper. Let $p$ be chosen as a fixed prime number for bosonic $p$-adic integral and as an odd prime number for fermionic $p$-adic integrals over $\mathbb{Z}_p$. The symbols $\mathbb{Z}_p, \mathbb{Q}_p$ and $\mathbb{C}_p$ represent the ring of $p$-adic rational integers, the field of $p$-adic rational numbers and the completion of an algebraic closure of $\mathbb{Q}_p$, respectively. Let $v_p$ denote the normalized exponential valuation of $\mathbb{C}_p$ with $|p|_p = \frac{1}{p}$; see [9,10,23] for more details.

Let $\mathrm{UD}(\mathbb{Z}_p)$ be the space of uniformly differentiable functions on $\mathbb{Z}_p$. For $f \in \mathrm{UD}(\mathbb{Z}_p)$, the $p$-adic $q$-integral (or $q$-Volkenborn integral) on $\mathbb{Z}_p$ is defined by Kim (see [9]) as follows:

$$I_q(f) = \int_{\mathbb{Z}_p} f(x) d\mu_q(x) = \lim_{N \to \infty} \frac{1}{[p^N]_q} \sum_{l=0}^{p^N-1} f(l) q^l \tag{1.15}$$

If we set $f(x) = e^{[x]_q v}$ in (1.15), it follows that

$$\int_{\mathbb{Z}_p} [x]_q^r d\mu_q(x) = \beta_{r,q}, \quad \text{see [9]}, \tag{1.16}$$

where $\beta_{r,q}$ represents Carlitz's $q$-Bernoulli numbers.

Let $C(\mathbb{Z}_p)$ be the space of all continous functions on $\mathbb{Z}_p$. For $f \in C(\mathbb{Z}_p)$, Kim introduced the fermionic $p$-adic $q$-integral on $\mathbb{Z}_p$ in [12] as follows:

$$I_{-q}(f) = \int_{\mathbb{Z}_p} f(x) d\mu_{-q}(x) = \lim_{N \to \infty} \frac{1}{[p^N]_{-q}} \sum_{l=0}^{p^N-1} f(l)(-q)^l. \tag{1.17}$$

By setting $f(x) = e^{[x]_q v}$ in (1.17), we see that

$$\int_{\mathbb{Z}_p} [x]_q^r d\mu_{-q}(x) = \mathcal{E}_{r,q}, \quad \text{see [12]}, \tag{1.18}$$

where $\mathcal{E}_{r,q}$ denotes the $q$-Euler numbers.

## 2. Main Results

In this section, we are now able to present the following theorems, which include probabilistic Bernoulli polynomials of higher order, probabilistic Stirling numbers of the second kind, probabilistic Euler polynomials, derivatives of probabilistic Bernstein polynomials and recurrence relation.

**Definition 2.1.** Let $B_{r,n}^Y(x, q)$ be the probabilistic $q$-Bernstein polynomials of degree $n$ given by the following generating function: For $x \in [0,1], n \in \mathbb{N}_0 := \mathbb{N} \cup \{0\}$ and $v \in \mathbb{C}$,

$$\sum_{n=r}^{\infty} B_{r,n}^{Y}(x,q) \frac{v^n}{n!} = \frac{\left(v[x]_q\right)^r}{r!} \left(\mathrm{E}\left[e^{vY}\right]\right)^{[1-x]_q}, \qquad \left(0 \le x \le 1; r \in \{0,1,2,\cdots,n\}\right).$$

**Remark 2.1.** In the limiting case as $q \to 1$ in Definition 2.1, it turns into Eq. (1.8) introduced by Karagenc *et al.* [8].

**Remark 2.2.** In the limiting case as $q \to 1$ and $Y = 1$, it reduces to the classical Bernstein polynomials.

**Theorem 2.1.** Let $Y$ be a random variable and $n \in \mathbb{N}_0$ with $r = \overline{0,n}$. The probabilistic $q$-Bernstein polynomials subject to a random variable $Y$ can be written as the probabilistic Bernoulli polynomials of higher order associated with $Y$:

$$B_{r,n}^{Y}(x,q) = \sum_{m=0}^{n} \binom{n}{m} \genfrac\{\}{0pt}{}{n-m}{r}_{Y} \beta_{m,Y}^{(r)}(1-[x]_{q^{-1}})[x]_q^r.$$

*Proof.* We first consider generating function of the probabilistic $q$-Bernstein polynomials associated with $Y$ in this form:

$$\sum_{n=r}^{\infty} B_{r,n}^{Y}(x:q) \frac{v^n}{n!} = [x]_q^r \frac{(\mathrm{E}[e^{Yv}]-1)^r}{r!} \left(\frac{v}{\mathrm{E}[e^{Yv}]-1}\right)^r (\mathrm{E}[e^{Yv}])^{(1-[x]_{q^{-1}})}.$$

Then,

$$= [x]_q^r \left(\sum_{n=r}^{\infty} \genfrac\{\}{0pt}{}{n}{r}_{Y} \frac{v^n}{n!}\right)\left(\sum_{n=0}^{\infty} \beta_{n,Y}^{(r)}((1-[x]_{q^{-1}})) \frac{v^n}{n!}\right).$$

It follows that

$$= \sum_{n=0}^{\infty} \left(\sum_{m=0}^{n} \binom{n}{m} \genfrac\{\}{0pt}{}{n-m}{r}_{Y} \beta_{m,Y}^{(r)}(1-[x]_{q^{-1}})[x]_q^r \right) \frac{v^n}{n!}$$

which yields the proofs when compared the coefficients $\frac{v^n}{n!}$.

Since

$$(\mathrm{E}[e^{Yv}])^{[x]_{q^{-1}}} = \sum_{r=0}^{\infty} \frac{(-[x]_{q^{-1}})_r}{r!} (\mathrm{E}[e^{Yv}]-1)^r,$$

we have the following theorem.

**Theorem 2.2.** Let $Y$ be a random variable and $n \in \mathbb{N}_0$ with $r = \overline{0,n}$. We have the explicit identity:

$$B_{r,n}^{Y}(x,q) = \sum_{l=0}^{n-r} \sum_{r=0}^{l} \binom{n}{r}\binom{n-r}{l} \genfrac\{\}{0pt}{}{l}{r}_{Y} (-[x]_{q^{-1}})_r [x]_q^r \mathrm{E}[Y^{n-r-l}].$$

*Proof.* From (1.9), we own

$$\sum_{n=r}^{\infty} B_{r,n}^Y(x,q) \frac{v^n}{n!} = \frac{(v[x]_q)^r}{r!} (E[e^{Yv}])^{1-[x]_{q^{-1}}}$$

$$= E[e^{Yv}] \frac{(v[x]_q)^r}{r!} (E[e^{Yv}])^{[x]_{q^{-1}}}$$

$$= E[e^{Yv}] \frac{(v[x]_q)^r}{r!} \sum_{r=0}^{\infty} (-[x]_{q^{-1}})_r \frac{(E[e^{Yv}] - 1)^r}{r!}$$

$$= E[e^{Yv}] \frac{(v[x]_q)^r}{r!} \sum_{r=0}^{\infty} (-[x]_{q^{-1}})_r \sum_{n=r}^{\infty} {n \brace r}_Y \frac{v^n}{n!}$$

$$= E[e^{Yv}] \frac{(v[x]_q)^r}{r!} \sum_{n=0}^{\infty} \left( \sum_{r=0}^{n} (-[x]_{q^{-1}})_r {n \brace r}_Y \right) \frac{v^n}{n!}$$

$$= \frac{(v[x]_q)^r}{r!} \left( \sum_{n=0}^{\infty} E[Y^n] \frac{v^n}{n!} \right) \left( \sum_{n=0}^{\infty} \left( \sum_{r=0}^{n} (-[x]_{q^{-1}})_r {n \brace r}_Y \right) \frac{v^n}{n!} \right)$$

$$= \frac{v^r [x]_q^r}{r!} \sum_{n=0}^{\infty} \left( \sum_{l=0}^{n} \left( \sum_{r=0}^{l} (-[x]_{q^{-1}})_r \binom{n}{l} {l \brace r}_Y E[Y^{n-l}] \right) \right) \frac{v^n}{n!}$$

$$= \sum_{n=r}^{\infty} \left( \sum_{l=0}^{n-r} \sum_{r=0}^{l} \binom{n}{r} \binom{n-r}{l} {l \brace r}_Y (-[x]_{q^{-1}})_r [x]_q^r E[Y^{n-r-l}] \right) \frac{v^n}{n!},$$

by matching the coefficients of $\frac{v^n}{n!}$, we arrive at

$$B_{r,n}^Y(x,q) = \sum_{l=0}^{n-r} \sum_{r=0}^{l} (-1)^r \binom{n}{r} \binom{n-r}{l} {l \brace r}_Y (-[x]_{q^{-1}})_r [x]_q^r E[Y^{n-r-l}],$$

which completes the proof.

**Theorem 2.3.** *Let $Y$ be a random variable, and $n \in \mathbb{N}_0$ with $r = \overline{0,n}$. Then we have*

$$[x]_q^r = \frac{1}{E[Y^{n-r}]} \sum_{l=0}^{n} \sum_{m=0}^{n-l} {n-l \brace m}_Y \frac{\binom{n}{l}}{\binom{n}{r}} ([x]_{q^{-1}})_m B_{r,l}^Y(x,q).$$

*Proof.* In the following alternative form of the generating function of probabilistic $q$-Bernstein polynomials:

$$\frac{(v[x]_q)^r}{r!} E[e^{Yv}] = (E[e^{Yv}])^{[x]_{q^{-1}}} \sum_{n=r}^{\infty} B_{r,n}^Y(x,q) \frac{v^n}{n!}.$$

On one hand,

$$I_1 = (E[e^{Yv}] - 1 + 1)^{[x]_{q^{-1}}} \sum_{n=0}^{\infty} B_{r,n}^Y(x,q) \frac{v^n}{n!}$$

$$= \sum_{m=0}^{\infty} ([x]_{q^{-1}})_m \frac{(E[e^{Yv}]-1)^m}{m!} \sum_{n=0}^{\infty} B^Y_{r,n}(x,q) \frac{v^n}{n!}$$

$$= \left(\sum_{m=0}^{\infty} ([x]_{q^{-1}})_m \frac{(E[e^{Yv}]-1)^m}{m!}\right)\left(\sum_{n=0}^{\infty} B^Y_{r,n}(x,q) \frac{v^n}{n!}\right)$$

$$= \left(\sum_{n=0}^{\infty} \sum_{m=0}^{n} ([x]_{q^{-1}})_m \left\{{n \atop m}\right\}_Y \frac{v^n}{n!}\right)\left(\sum_{n=0}^{\infty} B^Y_{r,n}(x,q) \frac{v^n}{n!}\right)$$

$$= \sum_{n=0}^{\infty} \left(\sum_{l=0}^{n} \sum_{m=0}^{n-l} \left\{{n-l \atop m}\right\}_Y \binom{n}{l} ([x]_{q^{-1}})_m B^Y_{r,l}(x,q)\right) \frac{v^n}{n!}.$$

On the other hand, it becomes

$$I_2 = \frac{(v[x]_q)^r}{r!} E[e^{Yv}]$$

$$= \sum_{n=r}^{\infty} \binom{n}{r} [x]_q^r E[Y^{n-r}] \frac{v^n}{n!}.$$

Since $I_1 = I_2$, we obtain the following equation by comparing coefficients $\frac{v^n}{n!}$ on both sides of the above:

$$[x]_q^r = \frac{1}{E[Y^{n-r}]} \sum_{l=0}^{n} \sum_{m=0}^{n-l} \left\{{n-l \atop m}\right\}_Y \frac{\binom{n}{l}}{\binom{n}{r}} ([x]_{q^{-1}})_m B^Y_{r,l}(x,q).$$

An immediate result of Theorem 2.3 when taken bosonic and fermionic $p$-adic integral both sides respectively of the above equations gives.

**Corollary 2.1** The following identity holds true:

$$\beta_{r,q} = \frac{1}{E[Y^{n-r}]} \sum_{l=0}^{n} \sum_{m=0}^{n-l} \left\{{n-l \atop m}\right\}_Y \frac{\binom{n}{l}}{\binom{n}{r}} \int_{\mathbb{Z}_p} ([x]_{q^{-1}})_m B^Y_{r,l}(x,q)\, d\mu_q(x),$$

$$\mathcal{E}_{r,q} = \frac{1}{E[Y^{n-r}]} \sum_{l=0}^{n} \sum_{m=0}^{n-l} \left\{{n-l \atop m}\right\}_Y \frac{\binom{n}{l}}{\binom{n}{r}} \int_{\mathbb{Z}_p} ([x]_{q^{-1}})_m B^Y_{r,l}(x,q)\, d\mu_{-q}(x).$$

**Theorem 2.4.** Let $Y$ be a random variable and $n \in \mathbb{N}_0, n \geq r$ with $r = \overline{0,n}$. Then we have

$$\sum_{j=0}^{n} \binom{n}{j} B^Y_{r,j}(x,q) \mathcal{E}^Y_{n-j}([x]_{q^{-1}}) = \binom{n}{r} [x]_q^r \mathcal{E}^Y_{n-r}(1).$$

*Proof.* Since

$$\frac{2}{E[e^{Yv}]+1} (E[e^{Yv}])^{[x]_{q^{-1}}} \sum_{n=r}^{\infty} B^Y_{r,n}(x,q) \frac{v^n}{n!} = \frac{v^r [x]_q^r}{r!} \frac{2E[e^{Yv}]}{E[e^{Yv}]+1}$$

$$= \left(\sum_{n=0}^{\infty} \mathcal{E}_n^Y([x]_{q^{-1}}) \frac{v^n}{n!}\right)\left(\sum_{n=0}^{\infty} B_{r,n}^Y(x,q) \frac{v^n}{n!}\right) = \sum_{n=0}^{\infty} \mathcal{E}_n^Y(1) \frac{v^{n+r}[x]_q^r}{r!\, n!}$$

$$= \sum_{n=0}^{\infty} \left(\sum_{j=0}^{n} \binom{n}{j} B_{r,j}^Y(x,q) \mathcal{E}_{n-j}^Y([x]_{q^{-1}})\right) \frac{v^n}{n!} = \sum_{n=r}^{\infty} \left(\binom{n}{r} [x]_q^r \mathcal{E}_{n-r}^Y(1)\right) \frac{v^n}{n!}.$$

Then, we possess

$$\sum_{j=0}^{n} \binom{n}{j} B_{r,j}^Y(x,q) \mathcal{E}_{n-j}^Y([x]_{q^{-1}}) = \binom{n}{r} [x]_q^r \mathcal{E}_{n-r}^Y(1).$$

**Theorem 2.5.** Let $Y$ be a random variable and $n \in \mathbb{N}_0$ with $r = \overline{0,n}$. Then we own

$$\sum_{j=0}^{n} \binom{n}{j} B_{r,j}^Y(x,q) \beta_{n-j}^Y([x]_{q^{-1}}) = \binom{n}{r} [x]_q^r \beta_{n-r}^Y(1).$$

*Proof.* Since

$$\frac{v}{\mathrm{E}[e^{Yv}] - 1} (\mathrm{E}[e^{Yv}])^{[x]_{q^{-1}}} \sum_{n=r}^{\infty} B_{r,n}^Y(x,q) \frac{v^n}{n!} = \frac{v^r [x]_q^r}{r!} \frac{v \mathrm{E}[e^{Yv}]}{\mathrm{E}[e^{Yv}] - 1}$$

$$= \left(\sum_{n=0}^{\infty} \beta_n^Y([x]_{q^{-1}}) \frac{v^n}{n!}\right)\left(\sum_{n=0}^{\infty} B_{r,n}^Y(x,q) \frac{v^n}{n!}\right) = \sum_{n=0}^{\infty} \beta_n^Y(1) \frac{v^{n+r}[x]_q^r}{r!\, n!}$$

$$= \sum_{n=0}^{\infty} \left(\sum_{j=0}^{n} \binom{n}{j} B_{r,j}^Y(x,q) \beta_{n-j}^Y([x]_{q^{-1}})\right) \frac{v^n}{n!} = \sum_{n=r}^{\infty} \left(\binom{n}{r} [x]_q^r \beta_{n-r}^Y(1)\right) \frac{v^n}{n!}.$$

Then, it becomes

$$\sum_{j=0}^{n} \binom{n}{j} B_{r,j}^Y(x,q) \beta_{n-j}^Y([x]_{q^{-1}}) = \binom{n}{r} [x]_q^r \beta_{n-r}^Y(1),$$

which completes the proof.

**Theorem 2.6.** The probabilistic $q$-Bernstein polynomials of degree $n$ can be expressed as the terms of the probabilistic $q$-Bernstein polynomials of degree $(n-1)$ as follows:

$$B_{r,n}^Y(x,q) = [x]_q B_{r-1,n-1}^Y(x,q) + [1-x]_q \mathrm{E}[Y] B_{r,n-1}^Y(x,q).$$

*Proof.* Taking the partial derivative with respect to $v$ on both sides gives

$$\sum_{n=r}^{\infty} B_{r,n}^Y(x,q) \frac{v^{n-1}}{(n-1)!} = \frac{[x]_q^r}{r!} \frac{\partial}{\partial v} \left(v^r (\mathrm{E}[e^{Yv}])^{[1-x]_q}\right)$$

$$= [x]_q \frac{(v[x]_q)^{r-1}}{(r-1)!} (\mathrm{E}[e^{Yv}])^{[1-x]_q} + [1-x]_q \mathrm{E}[Y] \frac{(v[x]_q)^r}{r!} (\mathrm{E}[e^{Yv}])^{[1-x]_q}$$

$$= \sum_{n=r-1}^{\infty} [x]_q B_{r-1,n}^Y(x,q) \frac{v^{n-1}}{(n-1)!} + \sum_{n=r}^{\infty} [1-x]_q \operatorname{E}[Y] B_{r,n}^Y(x,q) \frac{v^n}{n!}.$$

Thus, we arrive at the desired result.

Since

$$\log \operatorname{E}[e^{Yv}] = \sum_{n=0}^{\infty} \left( \sum_{l=0}^{n} (-1)^l \left\{ {n+1 \atop l+1} \right\}_Y \right) \frac{v^{n+1}}{(n+1)!},$$

we state the following theorem.

**Theorem 2.7.** Let $Y$ be a random variable and $n \in \mathbb{N}_0$ with $r = \overline{0,n}$. The partial derivative of the probabilistic $q$-Bernstein polynomials with respect to $x$ holds true:

$$\frac{\partial}{\partial x} B_{r,n}^Y(x,q) = nq^x \frac{\ln q}{q-1} B_{r-1,n-1}^Y(x,q) + q^{1-x} \frac{\ln q}{1-q} \sum_{j=0}^{n-1} \sum_{l=0}^{n-1-j} \binom{n}{j} B_{r,j}^Y(x,q) (-1)^l \left\{ {n-j \atop l+1} \right\}_Y.$$

*Proof.* In the generating function of the probabilistic $q$-Bernstein polynomials, when we employ the derivative operator with respect to $x$, it yields

$$\sum_{n=r}^{\infty} \frac{\partial}{\partial x} B_{r,n}^Y(x,q) \frac{v^n}{n!} = \frac{v^r}{r!} \frac{\partial}{\partial x} \left( [x]_q^r (\operatorname{E}[e^{Yv}])^{1-[x]_{q^{-1}}} \right)$$

$$= v \frac{(v[x]_q)^{r-1}}{(r-1)!} q^x \frac{\ln q}{q-1} (\operatorname{E}[e^{Yv}])^{(1-[x]_{q^{-1}})} + \frac{v^r}{r!} q^{1-x} \frac{\ln q}{1-q} (\operatorname{E}[e^{Yv}])^{(1-[x]_{q^{-1}})} \log \operatorname{E}[e^{Yv}]$$

$$= \sum_{n=r-1}^{\infty} q^x \frac{\ln q}{q-1} B_{r-1,n}^Y(x,q) \frac{v^{n+1}}{n!} + \sum_{n=r}^{\infty} q^{1-x} \frac{\ln q}{1-q} \log \operatorname{E}[e^{Yv}] B_{r,n}^Y(x,q) \frac{v^n}{n!}$$

$$= \sum_{n=r-1}^{\infty} q^x \frac{\ln q}{q-1} B_{r-1,n}^Y(x,q) \frac{v^{n+1}}{n!}$$

$$+ q^{1-x} \frac{\ln q}{1-q} \left( \sum_{n=r}^{\infty} B_{r,n}^Y(x,q) \frac{v^n}{n!} \right) \left( \sum_{n=0}^{\infty} \left( \sum_{l=0}^{n} (-1)^l \left\{ {n+1 \atop l+1} \right\}_Y \right) \frac{v^{n+1}}{(n+1)!} \right)$$

$$= \sum_{n=r-1}^{\infty} q^x \frac{\ln q}{q-1} B_{r-1,n}^Y(x,q) \frac{v^{n+1}}{n!}$$

$$+ q^{1-x} \frac{\ln q}{1-q} \sum_{n=0}^{\infty} \left( \sum_{j=0}^{n} \sum_{l=0}^{n-j} \binom{n+1}{j} B_{r,j}^Y(x,q) (-1)^l \left\{ {n-j+1 \atop l+1} \right\}_Y \right) \frac{v^{n+1}}{(n+1)!}.$$

By comparing coefficients $v^n$ on both sides of the above, we arrive at the desired result.

**Theorem 2.8.** Let $Y$ be a random variable and $n \in \mathbb{N}_0$ with $r = \overline{0,n}$.

$$F_{r;Y}^{(m)}(x,v,q) = \sum_{l=0}^{m} \binom{m}{l} [x]_q^l \frac{(r)_l}{r!} (r-l)! \operatorname{E}[Y^{m-l}][1-x]_q^{m-l} F_{r-l;Y}(x,v,q),$$

where $F_{r;Y}^{(m)}(x,v,q) = \frac{\partial^m}{\partial v^m} F_{r;Y}(x,v,q)$.

*Proof.* Using Leibniz's formula for the $m$th derivative with respect to $v$,

$$f(v) = \frac{(v[x]_q)^r}{r!} \text{ and } g(v) = (E[e^{Yv}])^{[1-x]_q}.$$

We obtain the following higher order partial derivative equation:

$$\frac{\partial^m}{\partial x^m} F_{r;Y}(x,v,q) = \sum_{l=0}^{m} \binom{m}{l} \left(\frac{\partial^l}{\partial x^l} f(v)\right) \left(\frac{\partial^{m-l}}{\partial x^{m-l}} g(v)\right),$$

substituting $f(v)$ and $g(v)$ into above, the theorem obtained.

## 3. Applications of Probabilistic $q$-Bernstein Polynomials in special random variables Y

In this section, employing specific random variables such as Poisson, Bernoulli, Binomial, Negative Binomial, Geometric, and Uniform distributions, we evaluate the probabilistic $q$-Bernstein polynomials. Additionally, we establish explicit identities, integral representations, and uncover new relationships among higher-order Bernoulli polynomials, Bell polynomials, higher-order Frobenius-Euler numbers, and Stirling numbers of the second kind.

**Theorem 3.1.** The probabilistic $q$-Bernstein polynomials associated with the Poisson random variable, $Y \sim Poisson(\alpha)$, can be written in terms of Bell polynomials as

$$B_{r,n}^Y(x,q) = \binom{n}{r} [x]_q^r \phi_{n-r}\big(\alpha(1-[x]_{q^{-1}})\big).$$

*Proof.* From (1.1), we have

$$\sum_{n=r}^{\infty} B_{r,n}^Y(x,q) \frac{v^n}{n!} = \frac{v^r [x]_q^r}{r!} e^{\alpha(1-[x]_{q^{-1}})(e^v-1)}$$

$$= \frac{v^r [x]_q^r}{r!} \sum_{n=0}^{\infty} \phi_n\big(\alpha(1-[x]_{q^{-1}})\big) \frac{v^n}{n!}$$

$$= \sum_{n=r}^{\infty} \binom{n}{r} [x]_q^r \phi_{n-r}\big(\alpha(1-[x]_{q^{-1}})\big) \frac{v^n}{n!},$$

by comparing the coefficients of $\frac{v^n}{n!}$ we arrive at the desired result.

Moreover, for the case of that $Y$ is a Poisson random variable, $Y \sim Poisson(\alpha)$, we own the following consequence which can be achieved through series manipulations.

**Theorem 3.2.** The following relation holds true:

$$B_{r,n}^Y(x,q) = \sum_{m=0}^{n-r} \binom{n}{r} \alpha^m [x]_q^r (1-[x]_{q^{-1}}) \begin{Bmatrix} n-r \\ m \end{Bmatrix}.$$

*Proof.* From (1.1), we have

$$\sum_{n=r}^{\infty} B_{r,n}^Y(x,q) \frac{v^n}{n!} = \frac{v^r [x]_q^r}{r!} e^{\alpha(1-[x]_{q^{-1}})(e^v-1)}$$

$$= \frac{v^r [x]_q^r}{r!} \sum_{m=0}^{\infty} \alpha^m (1-[x]_{q^{-1}}) \frac{(e^v-1)^m}{m!} = \sum_{m=0}^{\infty} \alpha^m (1-[x]_{q^{-1}}) \sum_{n=m}^{\infty} \left\{ \begin{matrix} n \\ m \end{matrix} \right\} \frac{v^n}{n!}$$

$$= \sum_{n=r}^{\infty} \left( \sum_{m=0}^{n-r} \binom{n}{r} \alpha^m [x]_q^r (1-[x]_{q^{-1}}) \left\{ \begin{matrix} n-r \\ m \end{matrix} \right\} \right) \frac{v^n}{n!}.$$

By comparing coefficients $\frac{v^n}{n!}$ on both sides of the above expression, we conclude the proof of this theorem.

From (1.7) and Theorem 3.2, we get the folllowing theorem.

**Corollary 3.1.** The following identity is true:

$$B_{r,n}^Y(x,q) = \sum_{m=0}^{n-r} \sum_{l=0}^{m} (-1)^l \alpha^m \left\{ \begin{matrix} n-r \\ m \end{matrix} \right\} \binom{n}{r} \binom{m}{l} q^{-mx} [x]_q^{r+l}.$$

**Corollary 3.2.** Taking $q$-Volkenborn integral in the Corollary 3.1, we have

$$\int_{\mathbb{Z}_p} B_{r,n}^Y(x,q) d\mu_q(x)$$

$$= \frac{\ln q}{(1-q)^{r+1}} \sum_{m=0}^{n-r} \sum_{l=0}^{m} \frac{1}{(1-q)^l} \sum_{j=0}^{l+r} (-1)^{l+j-1} \alpha^m \left\{ \begin{matrix} n-r \\ m \end{matrix} \right\} \binom{n}{r} \binom{m}{l} \binom{l+r}{j} \frac{j-m}{[j-m]_q}.$$

**Corollary 3.3.** Taking fermionic $p$-adic $q$-integral in the Corollary 3.1, we have

$$\int_{\mathbb{Z}_p} B_{r,n}^Y(x,q) d\mu_{-q}(x)$$

$$= \frac{2}{(1-q)^r} \sum_{m=0}^{n-r} \sum_{l=0}^{m} \frac{1}{(1-q)^l} \sum_{j=0}^{l+r} (-1)^{l+j} \alpha^m \left\{ \begin{matrix} n-r \\ m \end{matrix} \right\} \binom{n}{r} \binom{m}{l} \binom{l+r}{j} \frac{1}{1+q^{j-m}}.$$

Let $Y \sim Bernoulli(p_1)$ with the success probability $p_1$. Here is the following theorem.

**Theorem 3.3.** Let $Y$ be $Bernoulli(p_1)$. Then we have

$$B_{r,n}^Y(x,q) = [x]_q^r \binom{n}{r} \sum_{m=0}^{n-r} p_1^m (1-[x]_{q^{-1}})_m \left\{ \begin{matrix} n-r \\ m \end{matrix} \right\}.$$

*Proof.* From (1.2), we observe that

$$\sum_{n=r}^{\infty} B_{r,n}^Y(x,q) \frac{v^n}{n!} = \frac{v^r [x]_q^r}{r!} (p_1(e^v-1)+1)^{[1-x]_q}$$

$$= \frac{v^r[x]_q^r}{r!} \sum_{m=0}^{\infty} p_1^m \frac{(1-[x]_{q^{-1}})_m}{m!} (e^v - 1)^m$$

$$= \sum_{n=r}^{\infty} \left( [x]_q^r \binom{n}{r} \sum_{m=0}^{n-r} p_1^m (1-[x]_{q^{-1}})_m \begin{Bmatrix} n-r \\ m \end{Bmatrix} \right) \frac{v^n}{n!},$$

for the case of comparing coefficients $\frac{v^n}{n!}$ on the boths sides of the above, we thereby complete the proof of this theorem.

Let $Y \sim b(n, p_1)$ be binomial random variable with success $p_1$. Here is the subsequent theorem.

**Theorem 3.4.** The probabilistic $q$-Bernstein polynomials associated with the Binomial random variable $Y$ can be written in terms of linear combination of Stirling numbers of the second kind as

$$B_{r,n}^Y(x,q) = [x]_q^r \binom{n}{r} \sum_{m=0}^{n-r} p_1^m \left((1-[x]_{q^{-1}})(n-r)\right)_m \begin{Bmatrix} n-r \\ m \end{Bmatrix}.$$

*Proof.* It can be written from (1.3) that

$$\sum_{n=r}^{\infty} B_{k,n}^Y(x,q) \frac{v^n}{n!} = \frac{v^r[x]_q^r}{r!} (p_1(e^v - 1) + 1)^{n[1-x]_q}$$

$$= \frac{v^r[x]_q^r}{r!} \sum_{m=0}^{\infty} \binom{n(1-[x]_{q^{-1}})}{m} (p_1(e^v - 1))^m$$

$$= \frac{v^r[x]_q^r}{r!} \sum_{m=0}^{\infty} p_1^m \left(n(1-[x]_{q^{-1}})\right)_m \sum_{n=m}^{\infty} \begin{Bmatrix} n \\ m \end{Bmatrix} \frac{v^n}{n!}$$

$$= \sum_{n=r}^{\infty} \left( [x]_q^r \binom{n}{r} \sum_{m=0}^{n-r} p_1^m \left((1-[x]_{q^{-1}})(n-r)\right)_m \begin{Bmatrix} n-r \\ m \end{Bmatrix} \right) \frac{v^n}{n!},$$

which asserts the proof of this theorem.

Let $Y \sim Geo(p_1)$ be geometric distribution with obtaining the first success $p_1$.

**Theorem 3.5.** The probabilistic $q$-Bernstein polynomials associated with the Geometric random variable $Y$ can be written in terms of Frobenius-Euler number of higher order as

$$(-1)^{n-r} B_{r,n}^Y(x,q) = [x]_q^r \binom{n}{r} H_{n-r}^{(1-[x]_{q^{-1}})}(q_1).$$

*Proof.* Since

$$\sum_{n=r}^{\infty} B_{r,n}^Y(x,q) \frac{v^n}{n!} = \frac{v^r[x]_q^r}{r!} \left( \frac{p_1 e^v}{1 - q_1 e^v} \right)^{(1-[x]_{q^{-1}})}$$

$$= \frac{v^r [x]_q^r}{r!} \sum_{n=0}^{\infty} H_n^{(1-[x]_{q^{-1}})}(q_1) \frac{(-v)^n}{n!}$$

$$= \sum_{n=r}^{\infty} [x]_q^r \binom{n}{r} H_{n-r}^{(1-[x]_{q^{-1}})}(q_1) (-1)^{n-r} \frac{v^n}{n!},$$

we arrive at the desired result when compared the coefficients $\frac{v^n}{n!}$ on the above.

Using a Negative Binomial random variable $Y \sim NB(n, p_1)$, we derive the following theorem in a manner similar to that of the Geometric random variable.

**Theorem 3.6.** The probabilistic $q$-Bernstein polynomials associated with the Negative binomial random variable $Y$ can be expressed in terms of summation of the products of Frobenius-Euler number of higher order and Bernstein basis polynomials as:

$$B_{r,n}^Y(x,q) = \sum_{l=0}^{n} \binom{n}{l} a^{n-k} B_{r,l}(x) H_{n-l}^{(a(1-[x]_{q^{-1}}))}(q_1^{-1}).$$

Let $Y \sim U(0,1)$ be a uniform random variable. Then we finalize the paper with the following theorem.

**Theorem 3.7.** The probabilistic $q$-Bernstein polynomials associated with the uniform random variable $Y$ can be written in terms of Bernoulli polynomials of higher order as:

$$B_{r,n}^Y(x,q) = \binom{n}{r} [x]_q^r B_{n-k}^{([x]_{q^{-1}}-1)}.$$

*Proof.* By utilizing

$$\sum_{n=r}^{\infty} B_{r,n}^Y(x,q) \frac{v^n}{n!} = \frac{v^r [x]_q^r}{r!} \left(\frac{e^v - 1}{v}\right)^{[1-x]_q} = \frac{v^r [x]_q^r}{r!} \left(\frac{v}{e^v - 1}\right)^{[x]_{q^{-1}} - 1}$$

$$= \sum_{n=r}^{\infty} [x]_q^r \binom{n}{r} B_{n-r}^{([x]_{q^{-1}}-1)} \frac{v^n}{n!},$$

we achieve the desired result.

**Remark 3.1.** Each of these distributions Bernoulli, Binom, Poisson, Geometric, Negative Binomial, and Uniform has unique properties that make them suitable for use in generating functions. By utilizing these distributions, we can derive various types of formulas involving special polynomials. This approach proves useful in deriving new and interesting formulas in both mathematics and statistics.

## 4. Conclusion

Motivated by Adell and Lekuona ([2]) probabilistic viewpoint via the generating function of the probabilistic Stirling numbers of the second kind as

$$\sum_{n=m}^{\infty} S_Y(n,m) \frac{t^n}{n!} = \frac{\left(\mathrm{E}\left[e^{tY}\right]-1\right)^m}{m!},$$

several researchers have introduced different generalizations of special numbers and polynomials, see [3, 5, 7, 8, 9, 14, 15-21, 25-31]. In conclusion, we introduced a new class of polynomials, the probabilistic $q$-Bernstein polynomials, which extends the classical Bernstein polynomials by incorporating probabilistic and $q$-analysis. Through the generating function of these polynomials, we derived several important mathematical relations, linking probabilistic $q$-Bernstein polynomials with probabilistic Stirling numbers of the second kind, higher-order probabilistic Bernoulli polynomials associated with a random variable $Y$, Bell polynomials, Frobenius-Euler polynomials of higher order. We have also established recurrence and differentiation properties for probabilistic $q$-Bernstein polynomials. Using Leibniz's formula, an identity for the generating function was derived. In exploring various random variables, including the Poisson, Bernoulli, Binomial, Geometric, Negative Binomial, and Uniform distributions, we demonstrated how these choices lead to meaningful applications and relationships among probabilistic $q$-Bernstein polynomials, Bell polynomials, Stirling numbers of the second kind, higher-order Frobenius-Euler numbers, and higher-order Bernoulli polynomials. Additionally, by presenting $p$-adic $q$-integral and fermionic $p$-adic $q$-integral representations, we further extended the utility of probabilistic $q$-Bernstein polynomials within the context of $p$-adic analysis. Overall, this study not only introduces probabilistic $q$-Bernstein polynomials but also opens pathways for their application in analytic number theory, combinatorics, and mathematical analysis.

## 5. Future Recommendations

In the case $Y=1$, $\mathrm{E}\left[e^{vY}\right]$ equals exponential function, $e^v$. For this reason, when studying generating functions of special functions, polynomials and numbers, replacing $e^v$ with $\mathrm{E}\left[e^{vY}\right]$ introduced a new era of probabilistic extension. It is recommended that researchers explore probabilistic extensions of special functions and focus on deriving new and interesting relationships within subfields of mathematics and statistics.

**Acknowledgments**

The authors would like to thank to the anonymous reviewers and handling Editor for consecutive valuable comments and suggestions, which have improved the presentation of the paper substantially.

**Disclosure statement**

The authors declare no conflict of interest.

**Funding**

This research received no external funding.